\documentclass[12pt,a4paper,twoside]{article}
\usepackage{amsfonts}
\usepackage{amsthm,amsmath,amssymb}
\usepackage{color}
\def\trait #1 #2 #3 {\vrule width #1pt height #2pt depth #3pt}
\def\fin{
    \trait .3 5 0
    \trait 5 .3 0
    \kern-5pt
    \trait 5 5 -4.7
    \trait 0.3 5 0
\medskip}
\newtheorem{teor}{Theorem}[section]
\newtheorem{defin}[teor]{Definition}
\newtheorem{lemm}[teor]{Lemma}
\newtheorem{osse}[teor]{Remark}
\newtheorem{prop}[teor]{Proposition}
\newtheorem{defi}[teor]{Definition}
\newtheorem{coro}[teor]{Corollary}
\newtheorem{prob}[teor]{Problem}
\newtheorem{hypo}[teor]{Hypothesis}
\newcommand{\bele}{\begin{lemm}\begin{sl}}
\newcommand{\enle}{\end{sl}\end{lemm}}
\newcommand{\bedef}{\begin{defi}\begin{sl}}
\newcommand{\eddef}{\end{sl}\end{defi}}
\newcommand{\bete}{\begin{teor}\begin{sl}}
\newcommand{\ente}{\end{sl}\end{teor}}
\newcommand{\beos}{\begin{osse}\begin{rm}}
\newcommand{\eddos}{\end{rm}\end{osse}}
\newcommand{\bepr}{\begin{prop}\begin{sl}}
\newcommand{\empr}{\end{sl}\end{prop}}
\newcommand{\bepro}{\begin{prob}\begin{rm}}
\newcommand{\empro}{\end{rm}\end{prob}}
\newcommand{\bede}{\begin{defin}\begin{sl}}
\newcommand{\edde}{\end{sl}\end{defin}}
\newcommand{\beco}{\begin{coro}\begin{sl}}
\newcommand{\enco}{\end{sl}\end{coro}}
\newcommand{\behy}{\begin{hypo}\begin{sl}}
\newcommand{\enhy}{\end{sl}\end{hypo}}

\newcommand{\RR}{\mathbb{R}}

\newcommand{\beeq}[1]{\begin{equation}\label{#1}}
\newcommand{\eddeq}{\end{equation}}
\newcommand{\beeqa}[1]{\begin{eqnarray}\label{#1}}
\newcommand{\eddeqa}{\end{eqnarray}}
\newcommand{\beal}[1]{\begin{align}\label{#1}}
\newcommand{\eddal}{\end{align}}
\newcommand{\bespl}[1]{\begin{split}\label{#1}}
\newcommand{\edspl}{\end{split}}
\newcommand{\bega}[1]{\begin{gather}\label{#1}}
\newcommand{\edga}{\end{gather}}
\newcommand{\grad}{\nabla}
\newcommand{\beeqax}{\begin{eqnarray*}}
\newcommand{\eddeqax}{\end{eqnarray*}}
\def\qed{\ifmmode   \else \leavevmode\unskip\penalty9999 \hbox{}\nobreak\hfill
  \fi
  \quad\hbox{\hskip.5em\vrule width.4em height.6em depth.05em\hskip.1em}}
\def\endproofsym{\qed}

\def\endnobox{\def\endproofsym{}\end{proof}\def\endproofsym{\qed}}

\newcommand{\no}{\nonumber}
\newcommand{\beeqao}{\begin{eqnarray}\no}
\newcommand{\bealo}{\begin{align}\no}
\newcommand{\besplo}{\begin{split}\no}
\newcommand{\begao}{\begin{gather}\no}
\def\trait #1 #2 #3 {\vrule width #1pt height #2pt depth #3pt}
\def\fin{\hfill
    \trait .3 5 0
    \trait 5 .3 0
    \kern-5pt
    \trait 5 5 -4.7
    \trait 0.3 5 0
\medskip}

\newcommand{\duav}[1]{\langle{#1}\rangle}

\newcommand{\duavgamma}[1]{_{H^{-1/2}(\Gamma)}\langle{#1}\rangle_{H^{1/2}(\Gamma)}}
\newcommand{\duavomega}[1]{_{H^{-1}}\langle{#1}\rangle_{W^{1,2}_0}}

\newcommand{\vc}[1]{{\boldsymbol #1}}
\newcommand{\dt}{\partial_t}

\newcommand{\itt}{\int_0^t}

\newcommand{\io}{\int_\Omega}

\newcommand{\bd}{\boldsymbol{d}}

 \DeclareMathOperator{\dive}{div}

\let\TeXchi\chi
\def\chi{{\setbox0 \hbox{\mathsurround0pt
$\TeXchi$}\hbox{\raise\dp0 \copy0 }}}

\newcommand{\bu}{{\boldsymbol u}}

\newcommand{\Grad}{\nabla}
\newcommand{\tn}[1]{\mbox {\F #1}}

\font\F=msbm10   
scaled 1000

\newcommand{\ub}{{\boldsymbol u}}
\newcommand{\vb}{{\boldsymbol v}}

\def\fine{\hfill\kern4pt \vrule height4pt depth0pt width4pt }

\def\dive{\mbox{\rm div\,}}

scaled \magstep1 

   \numberwithin{equation}{section}
\begin{document}

\title{On a 3D isothermal model for nematic liquid crystals
accounting for stretching terms}

\author{Cecilia Cavaterra\thanks{The author was partially supported by MIUR-PRIN Grant 20089PWTPS
``Mathematical Analysis for inverse problems towards
applications''}\\
Dipartimento di Matematica, Universit\`{a} degli Studi di Milano\\
Via Saldini 50, 20133 Milano, Italy
\\
cecilia.cavaterra@unimi.it
\and
Elisabetta Rocca\thanks{The author was  supported by the FP7-IDEAS-ERC-StG
Grant \#256872 (EntroPhase)}\\
Dipartimento di Matematica, Universit\`{a} degli Studi di Milano\\
Via Saldini 50, 20133 Milano, Italy
\\
elisabetta.rocca@unimi.it}
\date{}
\maketitle

\begin{abstract}
\noindent In the present contribution we study a PDE system describing
the evolution of a nematic liquid crystals flow under kinematic transports
for molecules of different shapes.
More in particular, the evolution of the {\em velocity field} $\ub$ is
ruled by the Navier-Stokes incompressible system with a stress tensor
exhibiting a special coupling between the transport and the
induced terms. The dynamics of the {\em director field} $\bd$ is
described by a variation of a parabolic Ginzburg-Landau equation
with a suitable penalization of the physical constraint $|\bd|=1$.
Such equation accounts for both the kinematic transport by the
flow field and the internal relaxation due to the elastic energy.
The main aim of this contribution is to overcome the lack of a
maximum principle for the director equation and prove (without
any restriction on the data and on the physical constants of the
problem) the existence of global in time weak solutions under
physically meaningful boundary conditions on $\bd$ and $\ub$.
\end{abstract}
\medskip


\medskip

\noindent
{\bf AMS (MOS) subject clas\-si\-fi\-ca\-tion:}~~35D30, 35K45, 35Q30, 76A15.

\section{Introduction}
\label{sec:intro}

In this paper we consider a hydrodynamical system modeling the flow
of nematic liquid crystals.
Assuming that the material occupies a bounded spatial domain $\Omega
\subset \RR^3$ with a smooth boundary $\Gamma$,
the system couples the 3D incompressible Navier-Stokes
equations governing the motion of the velocities with a modified Allen-Cahn
equation for the director field, that is
\begin{align}
\label{incompressiINTRO}
&\dive\ub=0,  \hskip9.8truecm \text{in}\, (0, T) \times \Omega, \\
\label{momconsINTRO}
& \partial_t\ub+\dive(\ub\otimes\ub) + \Grad p =\dive\tn{T} + \vc{f},  \hskip5truecm \text{in}\, (0, T) \times \Omega, \\
\label{eqdINTRO}
&\partial_t \bd +\ub\cdot \Grad\bd-\alpha\bd\cdot\Grad\ub + (1-\alpha)\bd\cdot\Grad^T\ub=\gamma(\Delta\bd- \nabla_{\bd} W(\bd)),
\quad \text{in}\, (0, T) \times \Omega,
\end{align}
where
\begin{equation}\label{defTINTRO}
  \tn{T} = \tn{S} -\lambda\left(\Grad \bd \odot \Grad \bd\right)
  - \alpha\lambda(\Delta \bd-\grad_{\bd}W(\bd))\otimes\bd
  +(1-\alpha)\lambda\bd \otimes (\Delta \bd- \nabla_{\bd} W(\bd)),
\end{equation}
\begin{equation}\label{defSINTRO}
\tn{S} = \mu\left( \Grad \ub + \Grad^T \ub \right),
\end{equation}
$\tn{T}$ and $\tn{S}$ being the Cauchy stress and the Newtonian viscous stress tensors, respectively.
Here $\ub$ denotes the velocity field of the flow, $\bd$ is the director field and stands for the averaged
macroscopic/continuum molecular orientation in $\RR^3$, $p$ is a scalar function representing the hydrodynamic pressure (including
the hydrostatic part and the induced elastic part from the orientation field) and $\vc{f}$ is a given external force.
The positive constants $\mu$, $\lambda$ and $\gamma$ stand for the viscosity, the competition between kinetic energy and potential energy,
and the microscopic elastic relaxation time (Deborah number) for the molecular orientation field, respectively.
The function $W$ penalizes the deviation of the length $|\bd|$ from the value 1, which is due to liquid crystal molecules being
of similar size (cf. \cite{LinLiusimply}).  A typical example is a {\em double well potential}
as, e.g., $W(\bd)=(|\bd|^2-1)^2$. In general $W$ may  be written as a sum of a convex part, and a smooth,
but possibly non-convex one.
Finally, the constant $\alpha \in [0,1]$ is a parameter related to the shape of the
liquid crystal molecules. For instance, the spherical, rod-like and disc-like liquid crystal molecules correspond to the cases
$\alpha= \frac{1}{2}, \, 1$ and $0$, respectively (cf., e.g., \cite{Er}, \cite{Jeff} and \cite{sunliu}).

Concerning the notation, $\Grad_{\bd}$ represents the gradient with respect to the variable $\bd$.
$\Grad \bd \odot \Grad \bd$ denotes the $3\times 3$ matrix whose $(i,j)$-th entry is given by
$\nabla_i \bd \cdot \nabla_j \bd$, for $i\leq i,j\leq 3$, and $\otimes$ stands for the usual Kronecker product, i.e.,
$(\ub\otimes\ub)_{ij}:=\ub_i\ub_j$, for $i,j=1,2,3$. Finally, $\Grad^T$ indicates the transpose of the gradient.

We notice that this system was very successful in describing the coupling between the
velocity field $\bu$ and the director field $\bd$, especially in the liquid crystals of nematic type.

The hydrodynamics theory of liquid crystals was due to Ericksen and Leslie (cf. \cite{Er} and \cite{Le78}).
However, the general Ericksen-Leslie system was so complicated that only some special cases of it have been
investigated theoretically or numerically in the literature.

In this context, Lin and Liu (cf. \cite{LinLiusimply} and \cite{LL01}) formulated a simplified version of
the original model which has been analyzed also by several other authors, see, e.g., \cite{SW10}, \cite{S02},
\cite{wuxuliu}. In the simplified model, some meaningful physical terms, like the stretching and rotation
effects of the director field induced by the straining of the fluid, are not taken into account.

In a following paper by Coutand and Shkol\-ler \cite{CoutShkoller}, the authors considered a model in which the
stretching term is present and proved a local well-posedness result. Here, due to the presence of
the stretching term, the total energy balance does not hold.  To overcome such an inconvenience,
Sun and Liu \cite{sunliu} proposed a variant of the Lin and Liu model \cite{LinLiusimply} in which not only the
stretching term is included in the system, but also a suitable new component to the stress tensor is added.

In our paper, we refer to a slightly more general model derived by Wu, Xu and Liu in \cite{wuxuliu}.

The main interests in the topic come essentially from two directions.
First, the previous results in the literature were only obtained in 2D or in 3D, under
the assumption that the viscosity coefficient $\mu$ in the stress $\tn{S}$ (cf. \eqref{defSINTRO})
is sufficiently big with respect to proper norms of the initial data and with respect to other coefficients
like $\lambda$.
In the previous contributions \cite{sunliu} and \cite{wuxuliu} it was claimed that, due to
the impossibility of proving the boundedness in $L^\infty$ for $\bd$ (because the maximum principle cannot be
applied to the director field equation), the existence of solution was out of reach without assuming
to have a big viscosity coefficient in the velocity equation. Our main result shows that, even if
the existence of {\em classical solutions} cannot be proved  without any restriction on the size of
the coefficients and the data (as it is for the uncoupled 3D Navier-Stokes system), however, it is possible
to obtain the existence of {\em weak solutions}.
This is in agreement with the previous contributions in the field of incompressible 3D Navier-Stokes equations.

The main point here is an appropriate choice of the test functions leading to a
rigorous weak formulation of the system (cf. the following formulas  (\ref{weak1}--\ref{weak3})).
Let us notice that in the recent paper  \cite{cgr}  formal computations are performed in
order to show the existence of weak solutions for such a problem, but no rigorous
definition of the weak formulation, as well as no proof of existence of such solutions
are given.
In particular, in our manuscript, choosing properly the space of the test functions in the
weak momentum equation (cf. \eqref{weak2}), we obtain well-defined weak solutions.
This is necessary in order to deal with the stretching term in the Cauchy stress
$\tn{T}$ (cf. \eqref{defTINTRO}). More comments on this point are given in Remark~\ref{rem:weaksol}.

The second novelty of our analysis consists in the fact that, to our knowledge,  all the previous contributions in
the literature (except for \cite{cgr} where formal results are stated in case of Dirichlet, Neumann and periodic
boundary conditions) were obtained assuming periodic boundary conditions on the director field $\bd$.
However, from the applications point of view, the cases of non-homogeneous Dirichlet or Neumann boundary conditions
look more appropriate (cf., e.g., \cite{LS} where it is pointed out that the Neumann boundary conditions for $\bd$
are also suitable for the implementation of a numerical scheme).
Here we can rigorously deal with all the three types of conditions for $\bd$: periodic, Dirichlet and Neumann.

We note that another new aspect of this contribution relies on the techniques employed for the proof of
existence of solutions. This method is based on the combination of a Faedo-Galerkin approximation and a regularization
procedure, which is necessary in order to treat the high order stretching terms in the weak momentum equation.
Indeed, a non standard but also physically meaningful regularization of the momentum equation is obtained
by adding to it an $r$-Laplacian operator
acting on the velocities, i.e., we add in the stress tensor a term of the type $|\Grad \ub|^{r-2}\Grad \ub$.
The reader can refer to the series of papers on the  J.-L.  Lions models
(cf., e.g. \cite{lions1} and \cite[Chap. 2, Sec. 5]{lions2})
or on the Ladyzhenskaya models (cf., e.g., \cite{lad} and references therein),
where $|\Grad \ub|^r$ is replaced by $|\Grad \ub+\Grad^T \ub|^r$.

Let us mention that the global in time existence of weak solutions can be considered as a
starting point for the analysis of the long-time behavior of solutions. This problem, up to now, has been considered
for this system for instance in \cite{wuxuliu}, where the convergence of a global {\em strong solution} to a single steady
state as time tends to infinity has been proved in  2D  and in 3D for special sets of data, and in \cite{GrasWu}, where
the existence of a smooth global attractor of finite fractal dimension is obtained. However, in both these cases
the system was endowed with periodic boundary conditions.

Finally, our results have been recently used in \cite{prslong} and in \cite{ffrs}
where the authors prove, via \L ojasiewicz-Simon techniques, the
convergence of the trajectories to the stationary states of system
(\ref{incompressiINTRO}--\ref{eqdINTRO}) - coupled with suitable
boundary conditions - and the existence of weak solutions for a
non-isothermal system with Neumann (for $\bd$) and complete slip
(for $\ub$) boundary conditions, respectively.

\paragraph{Plan of the paper.} In the following Section~\ref{sec:model} we briefly introduce the modeling
approach leading to our system and we discuss the choice of the boundary conditions.
The {\em weak formulation} of (\ref{incompressiINTRO}--\ref{defSINTRO}) is given in Section~\ref{sec:mainres},
where the main theorem regarding existence of global in time solutions is stated. The proof is given in the two
remaining Sections~\ref{a} and \ref{A}. In particular, in Section~\ref{a}
the a priori estimates, from which we deduce a rigorous character of the approximated Faedo-Galerkin scheme
presented in Section~\ref{A}, are obtained.

\section{Mathematical model}
\label{sec:model}

In this section we briefly derive system (\ref{incompressiINTRO}--\ref{eqdINTRO}) from the macroscopic
point of view.
Hence, suppose that the material occupies a bounded spatial domain $\Omega
\subset \RR^3$, with a sufficiently regular boundary. Let
$\bu =\bu(t,x)$ denote the velocity in
the Eulerian reference system. Accordingly, the \emph{mass
conservation} is expressed by means of continuity equation
(i.e. the standard incompressibility constraint)
\begin{equation}\label{incompr}
\dive\ub=0 ,
\end{equation}
which is relevant in the context of nematic liquid crystals.

In the context of hydrodynamics, the basic variable is the flow map (the particle trajectory)
$x(X, t)$, where $X$ is the original labelling (the Lagrangian coordinate)
of the particle, also
referred to as the material coordinate, and $x$ is the current (Eulerian) coordinate, which is also called
the reference coordinate. For a given velocity field $\ub(x, t)$, the flow map is defined by the ODE and
initial condition
$$x_t=\ub(x(X, t), t),\quad x(X, 0)=X.$$
In order to incorporate the properties of the material, we need
to introduce the deformation tensor $\tn{F}$ such that
$$\tn{F}(X, t)=\partial_{X}(x(X, t)).$$
Combining these two equations  with the chain rule formula and defining, in Euler coordinates,
$\tilde{\tn{F}}(x, t)=\tn{F}(X, t)$, we obtain the following transport equation for $\tilde{\tn{F}}$
$$ \tilde{\tn{F}}_t+(\ub \cdot \Grad)\,\tilde{\tn{F}}=\Grad \ub \,\tilde{\tn{F}}.$$
Without ambiguity, in what follows we do not distinguish between
the two symbols $\tn{F}$ and $\tilde{\tn{F}}$. In this case the transport of the director field $\bd$ can be stated as
$$ \bd (x(X,t), t)=\tn{F}\, \bd_0(X),$$
where $\bd_0$ is the initial condition. Taking now the full time derivative of both sides
we get
$$\frac{d}{dt} \bd (x(X,t),t)=\frac{d\,\tn{F} }{dt}\,\bd_0(X)=\Grad \ub \,\tn{F}\,\bd_0(X)
=\Grad\ub\,\bd=(\bd\cdot\Grad)\,\ub.$$
This allows us to deduce that the {\em total transport} associated to the orientation field $\bd$ is
$$\partial_t \bd+\ub\cdot\Grad\bd-\bd\cdot\Grad\ub$$
and so we get  the following transport equation for $\bd$
\begin{equation}\label{direq}
\partial_t \bd +  \bu\cdot\Grad \bd-\alpha\bd\cdot \Grad\ub + (1 -\alpha)\bd\cdot \Grad^T\ub=\gamma(\Delta\bd-\nabla_{\bd} W(\bd)).
\end{equation}
In the general case, $W$ may be a penalty function that can be written as a sum of a convex (possibly non smooth) part, and a smooth,
but possibly non-convex one.
Equation (\ref{direq}) is associated with conservation of angular momentum. The left-hand side stands for the kinematic transport by the flow field,
while the right-hand side denotes the internal relaxation due to the elastic energy (cf., e.g., \cite{sunliu}).

Finally, following the lines of \cite{wuxuliu}, i.e., applying the Hamilton's principle to the
action functional
\[
{\cal A}(x)=\int_0^T\int_{U_0}\frac12|x_t|^2
-\lambda \left[\frac12|{\cal F}^{-T}\Grad{\cal F} \bd_0(X)|^2
+W({\cal F} \bd_0(X))\right]J\, {\rm d}X {\rm d} t,
\]
where $U_0$ is the region occupied by the fluid at time $t$ and $J={\rm det} \frac{\partial x}{\partial X}$,
we deduce the following conservation of linear momentum equation
\begin{align}\label{mombal}
&\partial_t\ub+\dive(\ub\otimes\ub) + \nabla p =\dive\Big(\mu\left( \Grad \ub + \Grad^T \ub \right)-\lambda\left(\Grad \bd \odot \Grad \bd\right) \\
\no
& -\alpha\lambda (\Delta \bd-\nabla_{\bd} W(\bd))\otimes\bd +(1-\alpha)\lambda \bd \otimes (\Delta \bd-\nabla_{\bd} W(\bd)) \Big) + \vc{f}.
\end{align}
This relation combines a usual equation describing the flow of an isotropic fluid with an extra nonlinear coupling term that is anisotropic.
The extra term is the induced elastic stress from the elastic energy through the transport, which is represented by the equation
for $\bd$.

On account of the previous analysis, we get system (\ref{incompressiINTRO}--\ref{eqdINTRO}). This system must be supplemented with a
suitable set of boundary conditions: the homogeneous Dirichlet boundary condition for the velocity field
\begin{equation}\label{slip}
\ub  = \mathbf{0} , \quad  \mbox{on} \ (0,T)\times\Gamma,
\end{equation}
together with the Neumann homogeneous boundary condition for the director field
\begin{equation}\label{neu}
\partial_{\vc{n}} \bd = \mathbf{0},  \quad \mbox{on} \ (0,T)\times\Gamma.
\end{equation}
Let us note that with similar techniques we can also treat the case of  Dirichlet
boundary condition for $\bd$
\begin{equation}\label{diri}
\bd|_{\Gamma} = \vc{h},  \quad \mbox{on} \ (0,T)\times\Gamma,
\end{equation}
assuming the boundary datum $\vc{h}$ regular enough.
In the following we will detail the proper modifications in the analysis in order
to treat also this case.

Finally, we can deal also with the case of periodic boundary conditions, but, since
in this case the computations are easier, we will not enter into full details.

\section{Main results}
\label{sec:mainres}

For the sake of simplicity, we restrict ourselves to the case
$\gamma=\lambda=1$.
Accordingly, our problem (\ref{incompressiINTRO}--\ref{eqdINTRO}),
endowed with initial and  boundary conditions, reads as follows
\begin{align}
\label{incompressiom}
&\dive\ub=0, \hskip9.7truecm \mbox{in} \ (0,T)\times\Omega, \\
\label{momconsom}
&\partial_t\ub+\dive(\ub\otimes\ub)+\nabla p=\dive(\mu\left( \Grad \ub + \Grad^T \ub \right)) -\dive \left( \Grad\bd\odot\Grad\bd \right ) \\
\nonumber
& -\dive{\left(\alpha(\Delta\bd-\nabla_{\bd} W(\bd))\otimes\bd - (1 - \alpha)\bd \otimes (\Delta\bd-\nabla_{\bd} W(\bd)) \right)}+ \vc{f}, \,
\mbox{in} \ (0,T)\times\Omega, \\
\label{eqdom}
&\partial_t \bd +\ub\cdot \Grad\bd - \alpha\bd\cdot\Grad\ub + (1 - \alpha)\bd\cdot\Grad^T\ub  =(\Delta\bd- \nabla_{\bd}W(\bd)),
\hskip0.6truecm \mbox{in} \ (0,T)\times\Omega, \\
\label{initial}
&\vc{u}(0, \cdot) = \vc{u}_0, \ \vc{d}(0, \cdot) = \vc{d}_0, \hskip6.9truecm \mbox{in} \ \Omega, \\
\label{slip1}
&\ub  = \mathbf{0}, \hskip10.3truecm  \mbox{on} \ (0,T)\times\Gamma, \\
\label{neu1}
&\partial_{\vc{n}} \bd = \mathbf{0},  \hskip9.9truecm \mbox{on} \ (0,T)\times\Gamma.
\end{align}
As already
mentioned in the previous section we will also consider the same
problem with Dirichlet boundary conditions for $\bd$, that is
\begin{equation}\label{dir1}
\bd_{|\Gamma}=\vc{h}, \quad \hbox{on }(0,T)\times \Gamma
\end{equation}
instead of condition \eqref{neu1}.

To begin with, we introduce a weak formulation of
(\ref{incompressiom}--\ref{neu1}) and state our main result on the
existence of global-in-time weak solutions, without any
restriction imposed on the initial data or on $\mu$.

\subsection{Weak formulation}
\label{weak}

In the weak formulation, the momentum equation (\ref{momconsom}) together
with the incompres\-sibility constraint (\ref{incompressiom}) are replaced by a family of integral identities
\begin{equation}
\label{weak1}
\int_{\Omega} \vc{u}(t, \cdot) \cdot \Grad \varphi = 0, \quad \mbox{for a.a.} \ t \in (0,T),
\end{equation}
\begin{equation} \label{weak2}
 \langle \partial_t\vc{u}, \varphi \rangle
- \int_{\Omega}\ub\otimes\ub :\Grad\varphi  + \io \mu\left( \Grad \ub + \Grad^T \ub \right):\Grad \varphi=
\end{equation}
\[
+ \int_{\Omega} (\Grad \bd\odot\Grad \bd) : \Grad \varphi + \alpha \io (\Delta\bd-\nabla_{\bd} W(\bd))\otimes\bd:\Grad\varphi
\]
\[
-(1- \alpha)\io \bd \otimes (\Delta\bd-\nabla_{\bd} W(\bd)):\Grad\varphi
+\int_\Omega {\vc f}\cdot \varphi, \quad \mbox{for a.a.} \ t \in (0,T), \]
for any
$\varphi \in W^{1,3}_0({\Omega}; \RR^3)$ such that $\dive\varphi=0$.

Equation (\ref{eqdom}) holds in the strong sense, thanks to the regularity obtained for $\bd$.
More specifically, we have
\begin{align}
\label{weak3}
&\partial_t \bd + \ub\cdot\Grad \bd -\alpha\bd\cdot\Grad\ub +(1 - \alpha)\bd\cdot\Grad^T\ub= \Delta\bd- \nabla_{\bd} W(\bd), \,
\mbox{ a.e. in }(0,T) \times \Omega,\\
\label{weak4}
&\partial_{\vc{n}}\bd = \mathbf{0},
\hskip8.6truecm \mbox{a.e. on }(0,T)\times\Gamma,  \\
\label{weak5}
&\bd(0, \cdot)=\bd_0, \hskip8truecm \mbox{a.e. in }\Omega.
\end{align}
A {\em weak solution} is a pair $(\ub,\,\bd)$ satisfying
(\ref{weak1}--\ref{weak5}), with $\ub(0,\cdot)=\ub_0$, a.e. in $\Omega$.

A {\em weak solution} of the Dirichlet problem is a pair $(\ub,\,\bd)$ satisfying
(\ref{weak1}--\ref{weak3}), (\ref{weak5}), with $\ub(0,\cdot)=\ub_0$, a.e. in $\Omega$, and
\begin{equation}\label{weak4dir}
\bd_{|\Gamma}=\vc{h}, \quad \hbox{a.e. on }(0,T)\times \Gamma\,.
\end{equation}

\subsection{Main existence theorems}

Before formulating the main result of this paper, let us state the
list of hypotheses imposed on the constitutive functions. We assume
that $\mu$ and $\alpha$ are positive coefficients, with $\alpha \in [0,1]$, and
\begin{align}
&W \in C^2(\RR^3), \quad W \geq 0,\label{hyp1}\\
&W=W_1+W_2 \, \hbox{ s.t. } \,  W_1 \hbox{ is convex and } W_2 \in C^1(\RR^3), \, \nabla W_2 \in C^{0,1}(\RR^3;\mathbb{R}^3)\label{hyp2}\\
&{\bf f} \in L^2(0,T; W^{-1, 2}(\Omega; \RR^3))\label{hyp3}.
\end{align}

Our first main result reads as follows.

\bete \label{theo1} Let $\Omega \subset \RR^3$ be a bounded domain
of class $C^{1,1}$. Assume that hypotheses (\ref{hyp1}--\ref{hyp3}) are  satisfied.
Finally, let the initial data be such that
\begin{align}
\label{hyp4}
&\vc{u}_0 \in W^{1,2}(\Omega; \RR^3),\quad  \dive \vc{u}_0 = 0 \ \mbox{in} \,  L^2(\Omega), \\
\label{hyp5}
&\vc{d}_0 \in W^{1,2} (\Omega; \RR^3), \quad W(\bd_0)\in L^1(\Omega).
\end{align}
Then problem (\ref{weak1}--\ref{weak5}) possesses a global in time  weak
solution ($\vc{u}$, $\vc{d}$)
belonging to the class
\begin{equation} \label{reg1}
\vc{u} \in L^\infty(0,T; L^2(\Omega; \RR^3))
\cap L^2(0,T;W^{1,2}_0(\Omega; \RR^3)),
\end{equation}
\begin{equation}\label{reg1bis}
\partial_t\vc{u} \in L^2(0,T;W^{-1,3/2}(\Omega; \RR^{3})),
\end{equation}
\begin{equation} \label{reg20}
W(\bd)\in L^\infty(0,T; L^1(\Omega)),\quad \Grad_{\bd} W(\bd) \in L^2((0,T)\times \Omega; \RR^3),
\end{equation}
\begin{equation} \label{reg2}
\vc{d} \in L^\infty(0,T;
W^{1,2}(\Omega; \RR^3))\cap L^2(0,T; W^{2,2}(\Omega;\RR^3))\cap H^1(0,T; L^{3/2}(\Omega;\RR^3)),
\end{equation}
and  additionally satisfying, for a.a.~$t\in(0,T)$, the\/
 {\rm energy inequality}
 \begin{align}\label{energy}
   &\frac{d}{dt} \io
   \big(| \bu |^2
    + |\nabla \bd|^2
    +2 W(\bd)\big)(t)
    + 2\big\| (- \Delta \bd+ \nabla_{\bd} W(\bd))(t) \big\|_{L^2(\Omega; \RR^3)}^2\\
    \no
    &+ {\mu} \| \nabla \bu(t) \|_{L^2(\Omega; \RR^{3\times3})}^2 \le C \|{\vc f}(t)\|^2_{W^{-1, 2}(\Omega; \RR^3)},
 \end{align}
 where $C$ denotes a positive constant depending on $\Omega$.
\ente

\beos\label{rem:weaksol}
 The regularity of the test function $\varphi$ in \eqref{weak2} can be justified by
 noting that, thanks to \eqref{reg1}, \eqref{reg20} and \eqref{reg2}, we have
 \[
\ub \in L^\infty(0,T; L^2(\Omega;\RR^3))\cap L^2(0,T;W^{1,2}_0(\Omega; \RR^3)),
\]
\[
\nabla \bd\in L^\infty(0,T; L^2(\Omega;\RR^{3\times 3}))\cap L^2(0,T;W^{1,2}(\Omega; \RR^{3\times 3})),
\]
\[
\Delta \bd -\Grad_{\bd} W(\bd) \in L^2((0,T)\times \Omega; \RR^3), \quad \bd \in L^\infty(0,T; W^{1,2}(\Omega;\RR^3)),
\]
and hence
 \begin{equation}\label{on:reg}
  \ub\otimes \ub,~~
   \nabla \bd\odot \nabla \bd,~~
    ( \Delta \bd - \nabla_{\bd} W(\bd)) \otimes \bd
   \in L^2(0,T;L^{3/2}(\Omega; \RR^{3\times 3})),
 \end{equation}
 whence their (distributional) divergence are elements
 of the space
 $$L^2(0,T;W^{-1,3/2}(\Omega; \RR^{3})).$$
 The same
 problem (cf. also \cite{sunliu} and \cite{wuxuliu}) occurs in
 the 2D case, where the same {\em weak formulation} is needed in order to get
 existence of {\em weak solutions}.
\eddos

Finally, we can state the second result in case of Dirichlet boundary conditions
for $\bd$.
\bete \label{theo2} Let $\Omega \subset \RR^3$ be a bounded domain
of class $C^{1,1}$. Assume that hypotheses (\ref{hyp1}--\ref{hyp2}) and (\ref{hyp4}--\ref{hyp5}) are  satisfied.
Suppose moreover that the assumptions
\begin{equation}\label{hyp6}
\vc{h}\in H^1(0,T;H^{-1/2}(\Gamma; \RR^3))\cap L^\infty(0,T; H^{3/2}(\Gamma;\RR^3)), \quad \vc{h}(0)={\bd_0}_{|\Gamma}
\end{equation}
hold true.
Then problem (\ref{weak1}--\ref{weak3}), (\ref{weak5}--\ref{weak4dir}) possesses a global in time  weak
solution ($\vc{u}$, $\vc{d}$)
belonging to the class stated in (\ref{reg1}--\ref{reg2}).
\ente

The rest of the paper is devoted to the proofs of Theorems
\ref{theo1} and \ref{theo2}.

In Section~\ref{a} we will prove the a-priori bounds on the solutions, detailing in Remark~\ref{rem:Dir}
the differences between the Neumann and the Dirichlet cases.
In Section~\ref{A} we introduce the regularization schemes and we perform the estimates on the approximated solutions in order to pass
to the limit in the approximated solutions. Since in the case of Dirichlet boundary conditions
\eqref{weak4dir} (instead of \eqref{weak4}) the technique is analogue, in Section~\ref{A} we perform the estimates and the approximation-passage to the limit
procedure on the Neumann system (\ref{weak1}--\ref{weak5}).

\section{A priori bounds}
\label{a}

We establish here a number of formal a priori estimates. These will assume a
rigorous character in the framework of the approximation scheme presented in
Section \ref{A} below. We inform the reader that a similar technique has been
used also in the subsequent work \cite{ffrs} (with respect to the present one)
for a non-isothermal model with different kind of boundary conditions.

Take $\varphi=\ub$ in \eqref{weak2} and test \eqref{weak3} by $-\Delta\bd+\nabla_{\bd} W(\bd)$ on $\Omega$.
Summing up the two resulting equalities, using the divergence theorem together with \eqref{weak1}, we obtain
\begin{equation} \label{a1}
\frac12\frac{\partial}{\partial t} \io\left(|\ub|^2 +|\Grad\bd|^2 + 2W(\bd)\right)
+ \mu\io|\Grad\ub|^2+\io|-\Delta\bd+\nabla_{\bd} W(\bd)|^2= \duavomega{ {\vc f},\ub} .
\end{equation}
Moreover, applying Schwarz and Poincar\'e inequalities on the right hand side, we can
deduce the energy estimate \eqref{energy}.
Integrating over
$(0,T)$ the inequality \eqref{a1}, and using
assumption \eqref{hyp3}, we get the a priori
bounds
\begin{equation}\label{apr1}
\vc{u} \in L^\infty(0,T; L^2(\Omega; \RR^3)) \cap L^2(0,T;
W^{1,2}(\Omega; \RR^3))\cap L^{10/3}((0,T)\times\Omega;\RR^3),
\end{equation}
\begin{equation}\label{apr2}
\vc{d} \in L^\infty (0,T; W^{1,2}(\Omega; \RR^3)),
\end{equation}
\begin{equation}\label{apr2bis}
-\Delta\vc{d}+\nabla_{\bd} W(\bd) \in L^2 (0,T; L^2(\Omega; \RR^3))\,.
\end{equation}

From \eqref{apr2bis}, on account of \eqref{neu} and \eqref{hyp1}, we get
\begin{equation*}
\int_0^T\int_\Omega |\Delta \vc{d}|^2 + \int_0^T\int_\Omega \Grad(\Grad_{\bd} W(\vc{d}))\nabla \vc{d}
=  \int_0^T\int_\Omega \vc{m} \cdot \Delta \vc{d},
\end{equation*}
$\vc{m}$ being a function in $L^2(0,T; L^2(\Omega; \RR^3))$.

Using once more assumption \eqref{hyp2}, from the previous equation we get
\begin{equation*}
\frac{1}{2} \int_0^T\int_\Omega |\Delta \vc{d}|^2 \leq \int_0^T\int_\Omega |\Grad(\Grad_{\bd}W_2(\vc{d}))||\nabla \vc{d}|
+ \frac{1}{2} \int_0^T\int_\Omega |\vc{m}|^2.
\end{equation*}
Recalling that $|\nabla \vc{d}| \in L^\infty (0,T; L^2(\Omega; \RR^3))$ (cf. \eqref{apr2}), then it holds (see again \eqref{apr2bis} and \eqref{hyp2})
\begin{equation}\label{apr3}
\bd\in L^2 (0,T; W^{2,2}(\Omega; \RR^3)), \quad  \nabla_{\bd} W(\bd)\in L^2((0,T)\times\Omega; \RR^3).
\end{equation}
From this result, since $\ub\cdot\Grad\bd$ and $\bd\cdot\Grad\ub$ belong to $L^2(0,T;L^{3/2}(\Omega; \RR^3))$, by comparison with \eqref{weak3}
we deduce
\begin{equation}\label{apr4}
\partial_t \bd\in L^2(0,T;L^{3/2}(\Omega;\RR^3)).
\end{equation}
Now, choosing $q(1-a)=2$ in the following  interpolation inequality
\begin{equation} \label{interp}
\|\Grad\bd\|_{L^{s}(\Omega; \RR^{3\times3})}^q\leq
c_1\|\Grad\bd\|_{L^2(\Omega; \RR^{3\times3})}^{aq}\|\Grad\bd\|_{L^6(\Omega;\RR^{3\times3})}^{(1-a)q}\,,
\end{equation}
holding true for
\begin{equation}\label{exp}
s,\,q\in [1,+\infty),\,\,a\in (0,1), \quad 1/s=(1-a)/6+a/2\,,
\end{equation}
and using (\ref{apr2}--\ref{apr3}), we get
\begin{equation}\label{apr4bis}
\Grad\vc{d} \in L^{4s/(3s-6)}(0,T;L^{s}(\Omega;\RR^{3\times3}))\,,
\end{equation}
which, taking $s=10/3$,  gives in particular
\begin{equation}\label{apr5}
\Grad\vc{d} \in L^{10/3}(0,T;L^{10/3}(\Omega;\RR^{3\times3})).
\end{equation}
This estimate turns out to be crucial for the proof of existence of solutions.
Hence, as a consequence of the previous estimates, we get
\begin{align}
\label{apr6}
&\left(-\left(\Grad \bd \odot \Grad \bd\right)
  +\alpha (\Delta \bd-\nabla_{\bd} W(\bd))\otimes\bd -
  (1-\alpha) \bd \otimes (\Delta \bd-\nabla_{\bd} W(\bd))\right)\\
\no
&\in L^{5/3}((0,T)\times\Omega;\RR^{3\times3})
\end{align}
and
\begin{align}
\label{apr7}
&\left(-\left(\Grad \bd \odot \Grad \bd\right)
  +\alpha(\Delta \bd-\nabla_{\bd} W(\bd))\otimes\bd - (1-\alpha) \bd \otimes (\Delta \bd-\nabla_{\bd} W(\bd))\right)\\
\no
&\in  L^2(0,T; L^{3/2}(\Omega; \RR^{3\times 3})).
\end{align}
Observe that from the a priori estimates \eqref{apr1}, \eqref{apr2} and \eqref{apr4} we derive the solution regularity
classes \eqref{reg1} and \eqref{reg2}, from which it follows \eqref{on:reg} and then \eqref{reg1bis}.
Moreover, we can prove that the set of the (weak) solutions to problem (\ref{weak1}--\ref{weak3})
is weakly stable (compact) with respect to these bounds, namely, taking any sequence of (weak) solutions satisfying
the above uniform bounds then it admits a convergent subsequence.
We omit the proof of the weak sequential stability, leaving the details to the reader, and
we devote the following section to the proof of Theorem \ref{theo1}. More precisely, we will construct a suitable family of
\emph{approximate} problems whose solutions weakly converge (up to subsequences)
to limit functions which solve the problem in the sense specified in Subsection~\ref{weak}.

\beos\label{rem:Dir} Let us detail here the different {\sl energy estimate} we obtain in the Dirichlet case.
If we we take non-homogeneous Dirichlet boundary conditions
for $\bd$, that is \eqref{weak4dir} instead of \eqref{weak4}, we get
\begin{equation}
\label{a2}
\frac12\frac{\partial}{\partial t} \io\left(|\ub|^2+|\Grad\bd|^2+2W(\bd)\right)
+\mu\io|\Grad\ub|^2+\io|-\Delta\bd+\nabla_{\bd} W(\bd)|^2
\end{equation}
\begin{equation}
\no
= \duavgamma{\vc{h}_t,\partial_{\vc{n}}\bd} + \duavomega{{\vc f},\ub}.
\end{equation}
Assuming  ${\vc h}$ satisfying \eqref{hyp6}, we can handle the first term on the right-hand side of \eqref{a2}
using standard trace theorems and
regularity results for elliptic equations (cf., e.g., \cite[Lemma 3.2, p. 263]{necas}) in this way
\begin{align}\label{a2bis}
\duavgamma{\vc{h}_t,\partial_{\vc{n}}\bd}&\leq C\|\vc{h}_t\|_{H^{-1/2}(\Gamma)}\|\bd\|_{H^2(\Omega)}\\
\no
&\leq C\left(\|\vc{h}_t\|^2_{H^{-1/2}(\Gamma)}+\|\vc{h}\|_{H^{3/2}(\Gamma)}^2\right)+\frac14\|\Delta\bd\|_{L^2(\Omega)}^2\,.
\end{align}
In order to treat the last term in \eqref{a2bis} we estimate from below the last integral in  \eqref{a2} as follows
\begin{align}\label{a2ter}
\|-\Delta\bd+\nabla_{\bd} W(\bd)\|_{L^2(\Omega)}^2&= \|\Delta \bd\|_{L^2(\Omega)}^2+\|\nabla_{\bd} W(\bd)\|_{L^2(\Omega)}^2
-2(\Delta \bd, \nabla_{\bd} W(\bd))\\
\no
& \geq \|\Delta \bd\|_{L^2(\Omega)}^2+\|\nabla_{\bd} W(\bd)\|_{L^2(\Omega)}^2 \\
\no
&\quad+ 2\io\Grad(\Grad_{\bd}W(\vc{d}))\nabla \vc{d} - 2\int_\Gamma \partial_n\bd(\nabla_{\bd} W(\bd))_{|\Gamma}\\
\no
&\geq \|\Delta \bd\|_{L^2(\Omega)}^2 - C\|\Grad \bd\|_{L^2(\Omega)}^2-\frac14\|\Delta\bd\|_{L^2(\Omega)}^2
\\
\no
&\quad-C\|\vc{h}\|_{H^{3/2}(\Gamma)}^2 - C\|\nabla_{\bd}W(\vc{h})\|_{L^2(\Gamma)}^2\,.
\end{align}
Here we have employed assumptions \eqref{hyp1}, \eqref{hyp2} and again standard elliptic estimates and trace theorems. Integrating over
$(0,T)$ the inequality \eqref{a2}, and using
assumptions \eqref{hyp1}, \eqref{hyp3}, \eqref{hyp6} together with (\ref{a2bis}--\ref{a2ter}) and a standard Gronwall lemma,
we get the a priori bounds
\begin{equation}\no
\vc{u} \in L^\infty(0,T; L^2(\Omega; \RR^3)) \cap L^2(0,T;
W^{1,2}(\Omega; \RR^3))\cap L^{10/3}((0,T)\times\Omega;\RR^3),
\end{equation}
\begin{equation}\no
\vc{d} \in L^\infty (0,T; W^{1,2}(\Omega; \RR^3)),
\end{equation}
\begin{equation}\no
-\Delta\vc{d}+\nabla_{\bd} W(\bd) \in L^2 (0,T; L^2(\Omega; \RR^3)),
\end{equation}
\begin{equation}
\no
\bd\in L^2(0,T; H^2(\Omega; \RR^3)), \quad \nabla_{\bd}W(\bd)\in L^2((0,T)\times\Omega; \RR^3)\,.
\end{equation}
Note that we have used, in particular, the  assumption $\vc{h}\in L^\infty(0,T; H^{3/2}(\Gamma;\RR^3))$ in order to get
$\vc{h}\in L^\infty(0,T; C^0(\overline\Gamma;\RR^3))$, which, thanks to the fact that $W\in C^2(\RR)$, implies
$\|\nabla_{\bd}W(\vc{h})\|_{L^2(\Gamma)}^2\in L^2(0,T)$. The remaining bounds (\ref{apr5}--\ref{apr7}) still holds true
in the Dirichlet case.
\eddos

\section{Approximations}

\label{A}

Here we introduce a double approximation scheme: a standard Faedo-Galerkin method is coupled with an
approximation  of the convective term and a regularization of the momentum equation by adding an $r$-Laplacian operator
acting on the velocities. Regarding the choice of these two regularization, let us note that for the first one (related to
the convective term) we follow the classical approach by Leray \cite{leray}, while the momentum equation with the $r$-Laplacian
acting on the velocity field $\ub$ (cf.  the following \eqref{approx1}) is a meaningful approximation of the standard Navier-Stokes system.
Indeed, one can refer to the series of papers on the  J.-L.  Lions models \cite{lions1} and \cite[Chap. 2, Sec. 5]{lions2}
or on the Ladyzhenskaya model (cf., e.g., \cite{lad} and references therein),
where $|\Grad \ub|^r$ is replaced by $|\Grad \ub+\Grad^T \ub|^r$.

The standard technique of the Faedo-Galerkin approximation scheme will be applied to construct the
solutions to the Navier-Stokes system (\ref{weak1}--\ref{weak2}) (see Temam~\cite{temam}).
To this aim, let us introduce first the Hilbert space
\[
W^{1,2}_{0,{\rm div}} = \{ \vc{v} \in W^{1,2}_0(\Omega; \RR^3) \ | \
\dive \vc{v} = 0, \ \mbox{a.e. in}\ \Omega \}
\]
and take an orthonormal basis $\{ \vc{v}_n \}_{n=1}^\infty$.
We fix $M, \,N \in \mathbb{N}$ such that $M\leq N$ and consider the finite-dimensional space
$X_N = {\rm span}\{ \vc{v}_n\}_{n=1}^N$.

Then the approximate velocity field $\vc{u}_{N,M} \in C^1([0,T]; X_N)$ solves
the Faedo-Galer\-kin system
\begin{equation} \label{approx1}
\frac{{\rm d}}{{\rm d}t} \int_\Omega \vc{u}_{N,M} \cdot \vc{v} =
\int_{\Omega} [ \vc{u}_{N,M} ]_M\otimes  \vc{u}_{N,M} : \Grad \vc{v}-\frac{1}{M}\io |\Grad\ub_{N,M}|^{r-2}\Grad\ub_{N,M}\cdot\Grad\vb
\end{equation}
\[
- \int_\Omega \mu \Big( \Grad \vc{u}_{N,M} + \Grad^T \vc{u}_{N,M} \Big): \Grad \vc{v}
+ \int_\Omega  \Grad \vc{d}_{N,M} \odot \Grad \vc{d}_{N,M} : \Grad \vc{v}
\]
\[
+ \alpha\int_\Omega\left(\Delta\bd_{N,M}-\nabla_{\bd} W(\bd_{N, M})\right)\otimes\bd_{N,M}: \Grad \vc{v}
\]
\[
- (1 - \alpha)\int_\Omega\bd_{N,M}\otimes \left(\Delta\bd_{N,M}-\nabla_{\bd} W(\bd_{N, M})\right): \Grad \vc{v}
+\int_\Omega \vc{f}\cdot\vc{v}, \quad \hbox{ for all }t \in [0,T],
\]
\begin{equation}
\label{approx1bis}
\int_\Omega \vc{u}_{N,M}(0, \cdot) \cdot \vc{v} = \int_\Omega
\vc{u}_0 \cdot \vc{v},
\end{equation}
for any $\vc{v} \in X_N$ and $r\in (3, 10/3)$. Here, the symbol $[ \vc{v} ]_M$ denotes
the orthogonal projection onto the space $X_M = {\rm span}\{\vc{v}_n\}_{n=1}^M$.
We observe that we need to introduce the additional term $\frac{1}{M}|\Grad\ub_{N,M}|^{r-2}\Grad\ub_{N,M}$ (cf. \eqref{weak2})
in order to obtain enough regularity for the velocity field in the director equation.

In \eqref{approx1} appears also the function $\vc{d}_{N,M}$ which is determined in terms of
$\vc{u}_{N,M}$ as the unique solution of the system
\begin{equation}
\label{approx2}
\partial_t \vc{d}_{N,M} + \vc{u}_{N,M} \cdot \Grad \vc{d}_{N,M} -\alpha \vc{d}_{N,M}\cdot \Grad\vc{u}_{N,M}
+(1 - \alpha) \vc{d}_{N,M}\cdot \Grad^T\vc{u}_{N,M}  +
\nabla_{\bd} W(\vc{d}_{N,M})
\end{equation}
\[
= \Delta \vc{d}_{N,M}, \quad\quad \hbox{ in $(0,T) \times \Omega$},
\]
\begin{equation}
\label{approx3}
\partial_\vc{n} \vc{d}_{N,M} = \mathbf{0},\quad \hbox{ on $(0,T) \times \Gamma$},
\end{equation}
\begin{equation} \label{approx4}
\vc{d}_{N,M}(0, \cdot) = \vc{d}_{0,M}, \quad \hbox{ in $\Omega$},
\end{equation}
$\vc{d}_{0,M}$ being a suitable smooth approximation of the initial datum $\vc{d}_0$ (cf. (\ref{weak3}--\ref{weak5})).

We will follow the original approach to the Navier-Stokes system by Leray \cite{leray} in order to regularize the convective
terms in (\ref{approx1}), (\ref{approx2}).
Hence, for any fixed $M,N$, we can solve problem (\ref{approx1}--\ref{approx4}) by means of a fixed point argument,
exactly as detailed in \cite[Chapter 3]{FN}.

Indeed, observe that all the a priori bounds derived formally in Section \ref{a}
still hold for our approximate problem. Hence, if we fix $\vc{u} \in C([0,T]; X_N)$, then we can find
$\vc{d} = \vc{d}[\vc{u}]$ solving (\ref{approx2}--\ref{approx4}). Inserting $\vc{d}[\vc{u}]$ in system
(\ref{approx1}--\ref{approx1bis}) we can define a mapping $\vc{u} \mapsto {\cal T}[\vc{u}]$, ${\cal T}[\vc{u}]$
being the solution of the system.
On account of the a priori bounds obtained in Section \ref{a}, we can easily show that ${\cal T}$ admits a fixed point
by means of the classical Schauder's argument on $(0,T_0)$, with $0< T_0 \leq T$.
Finally, applying again the a priori estimates, we are allowed to conclude that the approximate solutions can be
extended to the whole time interval $[0,T]$ (see \cite[Chapter 6]{FN} for details).

Now our strategy consists in passing to the limit first for $N\to\infty$ and then
for $M\to\infty$. This will be explained in the following subsections.


\subsection{Passage to the limit as $N\to\infty$}

On account of the regularizing term introduced in \eqref{approx1}, from the corresponding energy estimate we now obtain
\begin{equation}\label{unuova}
  M^{-1}\|\Grad\vc{u}_{N,M}\|^r_{L^r((0,T)\times\Omega;\mathbb{R}^{3\times 3})}\leq C, \ \mbox{ for } r \in (3, 10/3),
\end{equation}
from which we infer that, for any fixed $M$, the set of functions $|\Grad\vc{u}_{N,M}|^{r-2}\Grad\vc{u}_{N,M}$
is uniformly bounded in $L^{\frac{r}{r-1}}((0,T)\times\Omega;\mathbb{R}^{3\times3})$. Note that, since $r\in (3,10/3)$,
it holds $r/(r-1)\in (10/7, 3/2)$.

Hence, we deduce the following convergence results
\begin{align}
\label{cuN}
&\ub_{N,M}\to\ub_M \ \mbox{ weakly-(*) in } L^\infty(0,T;L^2(\Omega;\RR^3))
\cap L^2(0,T;W^{1,2}(\Omega;\RR^3)),\\
\label{cuNnuova}
&\Grad\ub_{N,M}\to\Grad\ub_M \ \mbox{ weakly in } L^{r}(0,T;L^r(\Omega;\RR^3)),\\
\label{cutN}
& \dt\ub_{N,M}\to\dt\ub_M \ \mbox{ weakly in }
  L^{\frac{r}{r-1}}(0,T; W^{-1, r/r-1}(\Omega;\RR^3)),\\
\label{cdN}
&\bd_{N,M}\to \bd_M \ \mbox{ weakly-(*) in } L^\infty(0,T;W^{1,2}(\Omega;\RR^3))
\cap L^2(0,T;W^{2,2}(\Omega;\RR^3)).
\end{align}
Observe that in \eqref{approx1} the projection on $X_M$ is kept in the convective term when passing
to the limit as $N \to \infty$.

Moreover, by virtue of \eqref{cdN} and a simple interpolation argument, we
have also
\begin{equation}\label{custg}
\Grad\bd_{N,M}\to \Grad \bd_M \ \mbox{ strongly
in }L^\eta((0,T)\times\Omega; \RR^{3\times3}),\ \mbox{ for }\eta\in [1,10/3).
\end{equation}
Going back to \eqref{cuN} and \eqref{cuNnuova}, by standard interpolation results,
some embedding properties of So\-bo\-lev spaces, and the Aubin-Lions lemma,
then we get
\begin{equation}\label{custrN}
  \ub_{N,M}\to\ub_M \ \mbox{ strongly in } L^s((0,T)\times\Omega;\RR^3), \ \mbox{ for some }  s >5.
\end{equation}
A combination of \eqref{custrN} with \eqref{custg} gives
\begin{equation}\label{convls}
 \vc{u}_{N,M}\cdot\Grad\bd_{N,M}\to \vc{u}_M\cdot\Grad \bd_M\
 \mbox{ strongly in } L^s((0,T)\times\Omega), \ \mbox{ for some }  s>2,
\end{equation}
whereas, combining  \eqref{cuNnuova} and \eqref{cdN}, one obtains
\[
\bd_{N,M}\cdot \Grad \ub_{N,M} \to \bd_M\cdot\Grad\ub_M \  \mbox{ weakly in } L^p((0,T)\times \Omega), \ \mbox{ for some } p>2.
\]
Moreover, it holds
\begin{equation}\label{cdtNbis}
 \dt \bd_{N,M}\to \dt \bd_M \ \mbox{ weakly in } L^2(0,T;L^{2}(\Omega;\RR^3)),
\end{equation}
and we have also
\[
|\Grad\ub_{N,M}|^{r-2}\Grad \ub_{N,M}\to \overline{|\Grad\ub_{M}|^{r-2}\Grad \ub_{M}}\
\mbox{weakly in } L^{r/r-1}((0,T)\times \Omega; \RR^{3\times3}).
\]
Considering now the pair of limit functions $(\vc{u}_M, \vc{d}_M)$, it can be proved that it solves the problem
\begin{equation}\label{aapprox2}
\int_{\Omega} \vc{u}_M(t, \cdot) \cdot \Grad \varphi = 0, \quad \mbox{for a.a.} \ t \in (0,T),
\end{equation}
\begin{equation} \label{aapprox1}
\int_0^t \duav{\partial_t\vc{u}_M, \varphi} -\int_0^t\int_{\Omega} \Big(  [\vc{u}_{M}]_M\otimes\vc{u}_M : \Grad \varphi \Big)
+\int_0^t \int_{\Omega} \mu \left(\Grad \ub_M + \Grad^T \ub_M\right): \Grad \varphi
\end{equation}
\[
= \int_0^t \int_{\Omega} \left( \Grad \bd_M \odot \Grad \bd_M + \alpha \left(\Delta\bd_M-\nabla_{\bd} W(\bd_M)\right)\otimes \bd_M\right) : \Grad \varphi
\]
\[
- \int_0^t \int_{\Omega}  (1 -\alpha) \bd_M \otimes \left (\Delta\bd_M-\nabla_{\bd} W(\bd_M)\right) : \Grad \varphi
\]
\[
- \frac1M \int_0^t\io \overline{| \Grad\vc{u}_{M} |^{r-2} \Grad\vc{u}_{M}}\cdot \nabla\varphi
+ \int_0^t \int_{\Omega} \vc{f}\cdot\varphi, \quad \mbox{for all} \ t \in (0,T),
\]
for any $\varphi \in C^\infty(\overline{\Omega}; \RR^3)$ such that $\dive \varphi  = 0$.

Passing to the limit as $N \to \infty$ also in the equation for $\vc{d}$, we get, a.e. in $(0,T)\times\Omega$,
\begin{equation} \label{aapprox3}
\partial_t \vc{d}_{M}
+ \vc{u}_{M} \cdot \Grad \vc{d}_{M} -\alpha \bd_M\cdot\Grad\ub_M
+(1 - \alpha) \bd_M\cdot\Grad^T\ub_M
 = \Delta \vc{d}_{M}-\nabla_{\bd}W(\vc{d}_{M}),
\end{equation}
as well as
\begin{equation} \label{aapprox4}
\partial_\vc{n} \vc{d}_{M} = \mathbf{0}, \quad \hbox{a.e.~in } (0,T) \times \Gamma,
\end{equation}
\begin{equation} \label{aapprox5}
\vc{d}_{M}(0, \cdot) = \vc{d}_{0,M}, \quad \hbox{a.e.~in } \Omega.
\end{equation}
Going back to \eqref{approx1}, if we replace $\vb$ by $\ub_{N,M}$ and then integrate with respect to time over $(0,t)$, we obtain
\begin{equation}\label{compaMN}
\|\ub_{N,M}(t)\|^2_{L^2(\Omega)}
+\itt\io\mu|\Grad \vc{u}_{N,M} + \Grad^T \vc{u}_{N,M}|^2
+\frac2M\itt\io |\Grad\vc{u}_{N,M}|^r
\end{equation}
\[
=\|\ub_0\|^2_{L^2(\Omega)} + 2\itt\io\ \left( \Grad \bd_{N,M} \odot \Grad \bd_{N,M}\right): \Grad \vc{u}_{N,M}
\]
\[
+ 2\alpha \int_0^t\io\left(\Delta\bd_{N,M}-\nabla_{\bd} W(\bd_{N,M})\right)\otimes \bd_{N,M}: \Grad \vc{u}_{N,M}
\]
\[
- 2(1 - \alpha) \int_0^t\io  \bd_{N,M} \otimes \left(\Delta\bd_{N,M}-\nabla_{\bd} W(\bd_{N,M})\right): \Grad \vc{u}_{N,M},
+ \int_0^t\io \vc{f}\cdot\vc{u}_{N,M},
\]
for all $t \in (0,T)$.

Next, on account of (\ref{cuN}--\ref{cutN}), we can
consider $\bu_M$ as a test function in \eqref{aapprox1} and get
\begin{equation}\label{compaM}
\|\ub_{M}(t)\|^2_{L^2(\Omega)}+\itt\io\mu|\Grad \vc{u}_{M} + \Grad^t
\vc{u}_{M}|^2+\frac2M\itt\io |\Grad\vc{u}_{M}|^r
\end{equation}
\[
=\|\ub_0\|^2_{L^2(\Omega)}+ 2\itt\io\ \left(\Delta\bd_{M}-\nabla_{\bd}W(\bd_{M})\right)\otimes \bd_{M} : \Grad \vc{u}_{M}
\]
\[
+ 2\alpha\int_0^t\io\left(\Delta\bd_{M}-\nabla_{\bd}W(\bd_{M})\right)\otimes \bd_{M} : \Grad \vc{u}_{M}
\]
\[
-2(1 -\alpha)\int_0^t\io  \bd_{M} \otimes \left(\Delta\bd_{M}-\nabla_{\bd} W(\bd_{M})\right): \Grad \vc{u}_{M}
+ 2 \int_0^t\io \vc{f}\cdot \vc{u}_{M}, \quad \mbox{for all} \ t \in (0,T).
\]
Actually, the $L^r$-regularity \eqref{cuNnuova} of $\Grad\bu_M$
is essential at this level since the terms
$\left(\Delta\bd_{M}-\nabla_{\bd}W(\bd_{M})\right)\otimes \bd_{M}$
and $\left(\Delta\bd_{M}-\nabla_{\bd}W(\bd_{M})\right)\otimes \bd_{M} $
do not necessarily belong to $L^2$. More precisely, the best we can conclude it is that
they lie in $L^{5/3}$ (cf.~\eqref{cdN} and \eqref{custg}).
Now, multiplying \eqref{approx2} by $\Delta\bd_{N, M} -\nabla_{\bd} W(\bd_{N,M})$ and integrating on $(0,t) \times \Omega$, we obtain
\begin{equation}
\label{compabisNM}
\|\Grad\bd_{N,M}(t)\|^2_{L^2(\Omega)} + 2\io W(\bd_{N,M})(t)+2\itt\io|\Delta\bd_{N,M}-\nabla_{\bd} W(\bd_{N,M})|^2
\end{equation}
\[
=\|\Grad\bd_{0,M}\|^2_{L^2(\Omega)}+2\io W(\bd_{0,M})+2\itt \left (\ub_{N,M}\cdot\Grad\bd_{N,M},\Delta\bd_{N,M}-\nabla_{\bd} W(\bd_{N,M})\right)
\]
\[
+2\itt\left(-\alpha \bd_{N,M}\cdot\Grad\ub_{N,M}
+(1 - \alpha) \bd_{N,M}\cdot\Grad^T\ub_{N,M},\Delta\bd_{N,M}-\nabla_{\bd} W(\bd_{N,M})\right),
\]
for all $t \in (0,T)$. If we multiply \eqref{aapprox3} by $\Delta \bd_{M} -\nabla_{\bd} W(\bd_M)$
and we integrate on $(0,t)\times\Omega$, then we deduce (for any $t\in (0,T)$)
\begin{equation}\label{compabisM}
\|\Grad\bd_{M}(t)\|^2_{L^2(\Omega)} + 2\io W(\bd_{M})(t)+2\itt\io|\Delta\bd_{M}-\nabla_{\bd}W(\bd_{M})|^2
\end{equation}
\[
=\|\Grad\bd_{0,M}\|^2_{L^2(\Omega)} + 2\io W(\bd_{0,M})
\]
\[
+ 2\itt \left ( \ub_{M}\cdot\Grad\bd_{M} - \alpha\bd_{M}\cdot\Grad\ub_{M}
+ (1 -\alpha)\bd_{M}\cdot\Grad^T\ub_{M} ,\Delta\bd_{M}-\nabla_{\bd} W(\bd_{M})\right ).
\]
Note that this is possible since \eqref{aapprox3} makes sense as a relation in $L^2$
(cf.~(\ref{convls}--\ref{cdtNbis})), thanks to the higher
regularity \eqref{cuNnuova} and \eqref{custrN} of $\bu_M$ and $\Grad\bu_M$ given by the regularizing
term $\frac{1}{M} |\Grad\vc{u}_{M}|^{r-2}\Grad\vc{u}_{M}$ in \eqref{approx1}.
Summing \eqref{compaMN} with \eqref{compabisNM} and then \eqref{compaM} with \eqref{compabisM},
passing to the limit as $N\to\infty$ we obtain
\[
\int_0^T\io |\Grad \ub_{N,M}|^r\to \int_0^T\io \overline{|\Grad\ub_{M}|^{r-2}\Grad\ub_M}:\Grad\ub_{M},
\]
\[
\int_0^T\io|\Delta\bd_{N,M}-\nabla_{\bd} W(\bd_{N,M})|^2\to \int_0^T\io|\Delta\bd_{M}-\nabla_{\bd}W(\bd_{M})|^2.
\]
So that, by means of standard Minty's trick and monotonicity argument, we infer
\[
\Grad\ub_{N,M}\to\Grad\ub_M\ \mbox{strongly in }L^{r}((0,T)\times\Omega;\RR^{3\times 3}),
\]
\[
\Delta \bd_{N,M}\to \Delta \bd_M\ \mbox{strongly in }L^{2}((0,T)\times\Omega;\RR^3).
\]
This concludes the passage to the limit as $N\to\infty$.


\subsection{Passage to the limit as $M\to\infty$}

The final step in the proof of the main result consists in passing to the limit as
$M\to\infty$ in (\ref{aapprox2}--\ref{aapprox5}).

First, we observe that we can still deduce the convergence results in (\ref{cuN}) and  (\ref{custrN}) when taking
$M\to \infty$.
Moreover, the following convergence results hold true
\begin{align}
\label{cdtuM}
& \dt\ub_M\to\dt\ub\ \mbox{ weakly in } L^{\frac{r}{r-1}}(0,T;W^{-1,r/r-1}(\Omega;\RR^3))\,,\\
\label{cdM}
&\bd_{M}\to \bd\ \mbox{ weakly-(*) in }  L^\infty(0,T;W^{1,2}(\Omega;\RR^3))\cap L^2(0,T;W^{2,2}(\Omega;\RR^3))\,,\\
\label{cdtM}
&\dt\bd_{M}\to \dt \bd \ \mbox{ weakly in } L^2(0,T;L^{3/2}(\Omega;\RR^3))\,,
\end{align}
and in particular
\begin{equation}
\label{cduMM}
M^{-1/(r - 1)}\Grad\ub_M\to 0 \mbox{ strongly in } L^{r - 1}((0,T)\times \Omega;\RR^{3\times3}).
\end{equation}
We are now in a position to pass to the limit as $M \to \infty$ in (\ref{aapprox2}--\ref{aapprox5}) and finally recover
(\ref{weak1}--\ref{weak5}).

\par
\medskip

\noindent
{\bf Acknowledgements.} We would like to thank Sergio Frigeri, Gianni Gilardi and the referees for their valuable comments and suggestions.

\par

\end{document}